\documentclass[aap,MSNbibl,seceqn,dvips]{arximspdf}

\doi{10.1214/13-AAP995} %kopijuoti is PTS
\volume{25}
\issue{1}
\pubyear{2015}
\firstpage{268}
\lastpage{286}

\makeatletter

\newtheorem{teo}{Theorem}[section]
\newtheorem{lem}[teo]{Lemma}
\newproclaim{ass}{Assumption}[section]
\newproclaim{rem}{Remark}[section]
\def\orig{o}
\newcommand{\ro}{\mathbf{r}}
\renewcommand{\i}{\mathbf{i}}
\newcommand{\n}{\mathbf{n}}
\newcommand{\GW}{\mathrm{GW}}
\newcommand{\UGW}{\mathrm{UGW}}
\newcommand{\BRW}{\mathrm{BRW}}
\makeatother

\begin{document}
\begin{frontmatter}

\title{Interacting growth processes and invariant~percolation}
\runtitle{Interacting growth processes and invariant percolation\hspace*{4pt}}

\begin{aug}
\author[A]{\fnms{Sebastian} \snm{M\"uller}\corref{}\ead[label=e1]{mueller@cmi.univ-mrs.fr}}
\runauthor{S. M\"uller}
\affiliation{Aix-Marseille Universit\'e (I2M)}
\address[A]{I2M, UMR 7373\\
CNRS\\
Aix-Marseille Universit\'e\\
Centrale Marseille\\
13453 Marseille\\
France\\
\printead{e1}} %adresu isvedimo komanda gale!
\end{aug}

% HISTORY:
\received{\smonth{4} \syear{2013}}
\revised{\smonth{7} \syear{2013}}

% ABSTRACT
%
\begin{abstract}
The aim of this paper is to underline the relation between reversible
growth processes
and invariant percolation. We present two models of interacting
branching random walks (BRWs), truncated BRWs and competing BRWs, where
survival of the growth process can be formulated as the existence of an
infinite cluster in an invariant percolation on a tree. Our approach is
fairly conceptual and allows generalizations to a wider set of
``reversible'' growth processes.
\end{abstract}

% KEYWORDS
\begin{keyword}[class=AMS]
\kwd[Primary ]{60J10}
\kwd{60J80}
\kwd[; secondary ]{05C80}
\end{keyword}
\begin{keyword}
\kwd{Survival}
\kwd{interacting branching random walk}
\kwd{invariant percolation}
\kwd{unimodular random networks}
\end{keyword}

\end{frontmatter}

%s1 #&#
\section{Introduction}\label{sec1}
We discuss two different interacting branching random walks (BRWs) in
discrete time. In the first model, called $\BRW_{N}$, only a finite
number $N$ of particles are allowed per site. A natural question is
whether the process $\BRW_{N}$ may survive with positive probability.
Partial answers to this question were given by Zucca \cite{Zucca11}.
We complete these results for symmetric BRWs on Cayley graphs in
Theorem~\ref{teotruncated}: $\BRW_{N}$ survives with positive
probability for sufficiently large $N$.

The second model concerns competing BRWs. Suppose there are two
different types or species of particles: invasive and noninvasive
particles. The invasive particles behave like particles in a usual BRW
and are not influenced by the noninvasive particles. These later,
however, die once they share a site with an invasive particle.
We prove (see Theorem~\ref{teojointsurvival1}) that for weakly
surviving (or transient) BRWs on Cayley graphs both processes may
coexist with positive probability.

Our proofs are based on a connection between the survival of reversible
growth processes and the existence of infinite clusters in percolation
on trees. This connection was used previously by Schramm \cite
{Kozma10} and Benjamini and Mueller \cite{benjamini11}. In the first
reference, BRWs are used to study connectivity properties of Bernoulli
percolation on nonamenable Cayley graphs. Benjamini and Mueller \cite
{benjamini11} use results on invariant percolation (on trees) to study
BRWs on unimodular random graphs.

In general, the study of interacting growth processes or particles
systems is challenging and a general treatment seems to be out of reach
(at least at the moment). The case-by-case study often involves a high
amount of technical effort. The approach given here is more conceptual
using \emph{soft proofs}. While we concentrate on two concrete
examples in this paper, we want to underline that our approach is
fairly general and only relies on two steps: the formulation of the
process as a unimodular random network and the control of the marginal
of the corresponding invariant percolation.

%s1.1 #&#
\subsection{Motivation}\label{sec1.1}
Besides highlighting the connection between growth processes and
percolation there are several other motivations for the underlying work.
One of these motivations is to propose models for spatial interaction
and competition of growth processes. One of the earliest and simplest
models of growth processes is the Galton--Watson branching process
where particles branch independently of the history of the process.
However, this may not be very realistic when there is competition for
limited resources such as space and food in the habitat. A considerable
effort was made to introduce dependence in the sense that the
individual reproduction may be influenced by the history of the
population. While many of these models consider dependence only on the
total population size (we refer to Kersting \cite{Kersting86} and the
paper referring to it for a mathematical introduction), the study of
models with local interactions is perhaps even more challenging. Models
in this direction are, for example, contact processes or restrained
BRWs. Here, particles breed depending on the local configurations of
the particles and one is interested in extinction, equilibrium and
explosion of the process. We refer to Bertacchi et al. \cite
{BPZ07} and references therein. A natural way to model local
dependencies between particles is also the truncated $\BRW_{N}$ that was
introduced in Bertacchi and Zucca \cite{BZ09} in continuous and in
Zucca \cite{Zucca11} in discrete time.%\footnote{The references given
%are far from being complete, but the abundance of existing literature
%does not allow to give a judicious presentation in this short note.}

As a byproduct of our approach, we can also control the following
processes. In a first model, dependencies are not only between
particles of the same generation but also between particles of
different generations. Suppose that each vertex has a finite amount of
resources that allow at most $N$ particles to branch; once the
resources are used any other particle visiting this site will die
without producing any offspring. For this model, an analogue result of
Theorem~\ref{teotruncated} holds for the weakly surviving regime; see
Section~\ref{sec3}. Another model of annihilating BRW that can be
treated is the process where the probability that two particles,
meeting at a same vertex, annihilate each other is a function of their
distance in the family tree. For instance, particles annihilate each
other if and only if their distance in the family tree is larger than
some constant $M$. Despite the nonmonotonicity of the model one can
prove that for $M$ large enough the process survives in the weakly
surviving regime.

Our second model describes two species competing for resources and
studies whether the weaker (or noninvasive) species has a chance of
survival. Models for competing spatial growth attracted a lot of
attention in the last decades. Perhaps the most common models in the
probability community are the voter model, the Richardson model and
mixtures of these two. We refer to H\"agstr\"om and Pemantle \cite
{HP98}, Kordzakhia and Lalley \cite{KL05} and Blair-Stahn \cite
{Blair-Stahn10} for an introduction and more references. Let us note
that most of the results are restricted to $\mathbb{Z}^{d}$ and make
strong use of a connection with first passage percolation.
Our model constitutes, to our knowledge, one of the first models beyond
$\mathbb{Z}^{d}$ and is more realistic for models where the space of
possible habitat grows exponentially (is expanding) in time. In
particular, this is relevant for models at the early stage of competing
populations. Moreover, it provides a stochastic model for the so-called
dominance displacement competition; we refer to Amarasekare \cite
{Amarasekare} for more details and references. In these kinds of
models, superior competitors can displace inferior competitors but not
vice versa. However, the inferior competitors can establish ``patches
or niches'' where the superior competitor does not colonize.
This latter phenomomen is highlighted by our theoretical result: as
long as the superior competitor does not colonize the whole space the
free patches are large enough to allow the inferior competitor to survive.

A somehow completely different motivation originates from a structure
theoretic question. Classification of groups in terms of the behavior
of random processes attracted a lot of attention. In particular, a
consequence of Gromov's famous theorem on groups of polynomial growth
is that a finitely-generated group admits a recurrent random walk if
and only if it contains a finite-index subgroup isomorphic to $\mathbb
{Z}$ or $\mathbb{Z}^{2}$, for example, see Chapter~3 in Woess
\cite{woess}. Kesten's criterion for amenability says that a
finitely-generated group is amenable if and only if the spectral radius
for all (or some) symmetric random walks is equal to $1$. This
phenomenon is also underlined by phase transitions on nonamenable
graphs whose study underwent a rapid development, for example,
see Lyons \cite{Lyons00}. Moreover,
Benjamini \cite{Benjamini02} proposed a deterministic competition
model that admits coexistence on hyperbolic groups but not on $\mathbb
{Z}^{d}$. A motivation for this is to find a stochastic process (or a
system of interacting processes) that shows an additional
phase-transition on (one-ended) hyperbolic groups compared to
nonhyperbolic groups. Theorem~\ref{teojointsurvival1} shows that
coexistence of competing BRWs, at least in the weakly surviving regime,
occurs regardless of the hyperbolicity of the underlying graph.

%s1.2 #&#
\subsection{Structure of the paper}\label{sec1.2}
We formulate our models and corresponding main results, Theorems \ref{teotruncated}~and~\ref{teojointsurvival1}, in the rest of
this section. In Section~\ref{secpreliminaries}, we introduce the
notation and basic results of random unimodular networks (URNs) and
present two preliminary results on percolation of URNs, Lemma \ref
{leminducedpercolation} and Theorem~\ref{teomarginal}. Section~\ref
{sec3} contains the proof of Theorem~\ref{teotruncated} and
Section~\ref{sec4} the one of Theorem~\ref{teojointsurvival1}.

%s1.3 #&#
\subsection{Branching random walks}\label{sec1.3}
The definition of a branching random walk (BRW) requires a probability
distribution, $\mu=(\mu_k)_{k\geq0}$, describing the branching and a
transition kernel, $P=(p(x,y))_{x,y\in V}$ describing the movement of
the particles on some underlying discrete space $V$. The BRW starts at
some initial position $\orig$ with one particle and then at each
(discrete) time step each particle splits into $k$ particles with
probability $\mu_{k}$ and each of the resulting particles moves one
step according to the transition kernel $P$. Both splitting and
movement of a particle at time $n$ are independent of the previous
history of the process and the behavior of the other particles at time $n$.

The expected number of offspring is denoted by $m=\sum_{k} k \mu_{k}$
and we will always assume that the corresponding Galton--Watson process
is supercritical, that is, $m>1$. Furthermore, we assume that
$P$ is the transition kernel of an irreducible random walk. Therefore,
the spectral radius $\rho=\rho(P)=\limsup_{n\to\infty}
(p^{(n)}(x,x))^{1/n}$, $x\in V$, of the underlying random walk is well defined.

There is an alternative description of BRWs that uses the concept of
tree-indexed random walks introduced in \cite{benjamini94}. Let
$(\mathbb{T},\ro)$ be a rooted infinite tree. Denote by $v$ the
vertices of $\mathbb{T}$ and let $|v|$ be the (graph) distance from
$v$ to the root~$\ro$. For any vertex $v$, denote $v^-$ the unique
predecessor of $v$, that is, $v^{-}\sim v$ and $|v^-|=|v|-1$. We denote
by $G=(V,E)$ a graph with vertex set $V$ and edge set $E$ and write
$(G,\orig)$ for a rooted graph.
%We often identify a Graph $G$ with its vertex set and denote the
%latter also by $G$.
The Cayley graph of a finitely generated group $\Gamma$ with respect
to some generating set $S$ will also be denoted by $G$; in this case,
$V=\Gamma$.

Let $(G,\orig)$ be a rooted graph. The tree-indexed process
$(S_v)_{v\in\mathbb{T}}$ is defined inductively such that $S_\ro
=\orig$ and for vertices $x,y\in V$ we have
\begin{eqnarray*}
\mathbb{P} \bigl(S_v=x| S_{v^-}=y, \bigl
\{S_{w}\dvtx w\notin\bigl\{v, v^{-} \bigr\}, |w|\leq n \bigr\}
\bigr) &=&\mathbb{P}(S_v=x| S_{v^-}=y)
\\
&=&p(x,y).
\end{eqnarray*}

A tree-indexed random walk becomes a BRW if the underlying tree
$\mathbb{T}$ is a realization of a Galton--Watson process. We call
$\mathbb{T}$ the family tree and $G$ the base graph of the BRW. The
symbol $\mathbb{T}$ will sometimes stand for random variables taking
values in the space of trees and sometimes for their realizations.

If $G$ is a Cayley graph, the BRW can be described as a marked (or
labeled) Galton--Watson tree. Let $\Gamma$ be a finitely generated
group with group identity $\orig$ and write the group operations
multiplicatively. Let $q$ be a symmetric probability measure on a
finite symmetric generating set of $\Gamma$. The corresponding random
walk on $\Gamma$ is the Markov chain with state space $\Gamma=V$ and
transition probabilities $p(x,y)=q(x^{-1}y)$ for $x,y\in\Gamma$.
Equivalently, the random walk (starting in $x$) can be described as
\[
S_n=x X_1\cdots X_n,\qquad n\geq1,
\]
where the $X_i$ are i.i.d. random variables with distribution $q$.
In order to define the BRW, we label the edges of $\mathbb{T}$ with
i.i.d. random variables $X_v$ with distribution $q$; the random
variable $X_v$ is the label of the edge $(v^-,v)$. These labels
correspond to the steps of the BRW and positions of particles are given
by $S_v=\orig\cdot\prod_{i} X_{v_i}$ where $\langle v_0=\ro, v_1,
\ldots, v_n=v\rangle$ is the unique geodesic from $\ro$ to $v$ at
level $n$.

A BRW is said to \textit{survive strongly} (\textit{or locally}) if
every vertex is
visited infinitely many times with positive probability and to \textit
{survive weakly} if the process survives with positive probability but
every finite subset is eventually free of particles. In formul\ae,
\begin{eqnarray*}
\mbox{strong survival} &\Leftrightarrow&\forall x\in G\dvtx \mathbb{P}
\bigl(\bigl|\{ v\dvtx
S_{v}=x\}\bigr|=\infty\bigr)>0,
\\
\mbox{weak survival} &\Leftrightarrow&\mathbb{P}\bigl(|\mathbb{T}|=\infty
\bigr)>0\quad\mbox{and}\quad \forall x\in G\dvtx \mathbb{P} \bigl(\{v\dvtx S_{v}=x\}|=\infty
\bigr)=0.
\end{eqnarray*}
Important to note that several authors speak sometimes of \textit
{transience} and
\textit{recurrence} of BRWs instead of weak and strong
survival. A consequence of the classification of recurrent groups and
Kesten's amenability criterion is that a BRW on a Cayley graph survives
strongly if and only if $m\rho(P)> 1$; see also \cite{gantert04} for
an alternative proof.

We make the following standing assumptions on the underlying
probability measures.

%as1.1 #&#
\begin{ass}\label{ass}
\begin{itemize}
\item The underlying Galton--Watson process is supercritical, $m=\sum
_{k} k \mu_{k}>1$, and the offspring\vspace*{1.5pt} distribution $\mu$ is of finite
support, that is, there exists some $d$ such that $\sum_{k=0}^{d-1} \mu
_{k}=1$. Furthermore, we assume that $\mu_{1}>0$.
\item Let $G$ be a finitely generated group with symmetric finite
generating set $S$. The distribution $q$ of the random walk on $G$ is
symmetric and such that $\operatorname{supp}(q)=S$ and $q(e)>0$.
\end{itemize}
\end{ass}

%re1.1 #&#
\begin{rem}
While the assumptions on supercriticality of the Galton--Watson process
and symmetry and irreducibility of the random walk are crucial, the
other assumptions are made for sake of a better presentation and to
avoid periodic subtleties.
In particular, the assumption that the genealogy of the invasive BRW is
a Galton--Watson process with finite support is not needed anywhere.
Moreover, the assumption on finite support of the offspring
distribution can be removed from the $\BRW_{N}$ and the noninvasive
process by adding an additional coupling of the Galton--Watson process
with a ``truncated'' version.
For instance, denote by $\mu^{(M)}$ a truncated version of $\mu$,
that is, $\mu^{(M)}_{k}=\mu_{k}$ for all $k< M$ and $\mu
^{(M)}_{M}=\sum_{k=M}^{\infty} \mu_{k}$, where $M$ is chosen
sufficiently large such that $m^{(M)}>1$.
Theorem~\ref{teotruncated} guarantees the existence of some
$N_{c}^{(M)}$ such that the $\BRW_{N}^{(M)}$ with underlying offspring
distribution $\mu^{(M)}$ may survive if $N\geq N_{c}^{(M)}$. Since
$\BRW_{N}^{(M)}$ is stochastically dominated by $\BRW_{N}$ this implies
survival of the latter if $N\geq N_{c}^{(M)}$. The argument for the
offspring distribution of the noninvasive process is similar.
\end{rem}

%s1.4 #&#
\subsection{Truncated branching random walk}\label{sec1.4}
Branching random walks may be used to describe the evolution of a
population or particle system at early stage with no restrictions on
resources. In order to refine the model, one might introduce a limit of
particle at each site: for some $N\in\mathbb{N}$ at most $N$
particles are allowed at a same site at the same time. While most of
the existing models describing variants of this models are in
continuous time, we prefer a description in discrete time since our
proof technique is more suitable to this setting. Let $N\in\mathbb
{N}$ and $W$ be a finite set and denote by $C(W,N)$ a random variable
that chooses uniformly $N$ elements from the set $W$ with the
convention that if $N>|W|$ then all $|W|$ elements are chosen. %All
%random variables of the kind $C(V,N)$ that will appear are supposed to
%be independent of everything else.

We define an auxiliary process: let $S_{v}^{\mathrm{aux}}$ be a BRW with
offspring distribution $\mu$ and transition kernel $P$ and denote by
$\mathbb{T}^{\mathrm{aux}}$ the corresponding family tree. For every $x\in X$
and $n\in\mathbb{N}$, we denote by $W_{n,x}=\{w\dvtx |w|=n,
S_{w}^{\mathrm{aux}}=x\}$ the particles at generation $n$ that are in position
$x$. We add a special state $\dagger$ to the state space $V$ and
define the process $\BRW_{N}$ on $V\cup\{ \dagger\}$ as
%
%e1.1 #&#
\begin{equation}
S_{v}^{N}=\cases{ S_{v}^{\mathrm{aux}}, &\quad
if $v\in C(W_{|v|,S_{v}^{\mathrm{aux}}},N)$,
\cr
\dagger, &\quad otherwise,}\qquad\mbox{for
all }v\in\mathbb{T}^{\mathrm{aux}},
\end{equation}
where $\{C(W_{n,x},N), n\in\mathbb{N}, x\in V\}$ is a family of
independent random variables independent of $S_{v}^{\mathrm{aux}}$.
The state $\dagger$ induces a site percolation on the family tree
$\mathbb{T}^{\mathrm{aux}}$ in the following way: declare a vertex $v$ closed
if $S_{v}^{\mathrm{aux}}=\dagger$ and open otherwise. Configurations of this
percolation are denoted by $\eta_{\dagger}$, where $\eta_{\dagger
}(v)=1$ corresponds to the fact that the site $v$ is open and $\eta
_{\dagger}(v)=0$ to the fact that the site $v$ is closed.
We denote by $\mathbb{T}_{\ro}$ the connected component containing
the root.

The truncated process $\BRW_{N}$ can therefore be denoted by
$(S^{N}_{v})_{v\in\mathbb{T}_{r}}$. We say that $\BRW_{N}$ \textit
{survives} if $|\mathbb{T}_{\ro}|=\infty$.
It is easy to see that survival is a monotone property: the process
survives for $N_{2}>N_{1}$ if it survives for $N_{1}$. The following
results asserts that there exists a nontrivial phase transition.

%th1.1 #&#
\begin{teo}\label{teotruncated}
Let $P$ be the transition kernel of a symmetric irreducible random walk
on an infinite finitely generated group $\Gamma$ and let $\mu$ be an
offspring distribution with $m>1$. Then there exists a critical value
$N_{c}=N_{c}(\mu,P)<\infty$ such that if $N\leq N_{c}$ the process
dies out a.s. and if $N>N_{c}$ the process survives with positive probability.
\end{teo}

%re1.2 #&#
\begin{rem}
Some of the results have been proved by Zucca \cite{Zucca11} in
Theorem 6.5: the case $m>1/\rho(P)$ was settled completely but the
case $m\leq1/\rho(P)$ was only proved for some BRWs on $\mathbb
{Z}^{d}$ and on the homogeneous tree. Zucca's results are presented in
the more general context of quasi-transitivity and treat some BRWs with
drift on $\mathbb{Z}^{d}$ that are not covered by our result. While
his proof technique is different to ours, it is interesting to note
that he uses as well a percolation argument for the case $m>1/\rho
(P)$; this time directed percolation on products of~$\mathbb{N}$.
\end{rem}

%re1.3 #&#
\begin{rem}
The notion of weak and strong survival can be adapted to the truncated
BRW in a natural way. If the underlying BRW survives weakly then
$\BRW_{N}$ survives weakly if $N>N_{c}$. If the underlying BRW survives
strongly then Theorem 6.5 in \cite{Zucca11} implies that there exists
some $N_{c}^{(s)}$ such that $\BRW_{N}$ survives strongly if
$N>N_{c}^{(s)}$. However, it is not known if there is an additional
regime of weak survival in this case, that is, there exists one
$N_{c}^{(w)}<N_{c}^{(s)}$ such that $\BRW_{N}$ survives weakly if
$N_{c}^{(w)}< N \leq N_{c}^{(s)}$.
\end{rem}

%s1.5 #&#
\subsection{Competing branching random walks}\label{sec1.5}
We consider two competing BRWs that interact in the following way. One
BRW is \textbf{i}nvasive, that is, the particles are not
influenced by the other particles, and the second is \textbf
{n}oninvasive in the sense that once a particle shares a site (at the
same time) with an invasive particle it dies without having any
offspring. The particles live on an infinite finitely generated group
$\Gamma$ and we note $(\mu_{\i},P_{\i}), (\mu_{\n},P_{\n})$ for
their offspring distribution and transition kernels. Moreover, denote
by $m_{\i}$ and $m_{\n}$ their expected number of offspring.

We give a formal definition of a slightly different process in the following.
The branching distributions $\mu_{\i}$ and $\mu_{\n}$ give rise to
two family trees $\mathbb{T}^{\i}$ and $\mathbb{T}^{\n}$. The
noninvasive BRW will start in $\orig$ and the invasive in some point
$x\neq\orig$.

The invasive BRW $(S_{v}^{\i})_{v\in\mathbb{T}^{\i}}$ is defined as
an ordinary BRW. In order to define the noninvasive BRW, we first
construct an intermediate version. Let $S_{v}^{\mathrm{aux}}$ be an ordinary BRW
with $(\mu_{\n},P_{\n})$ and introduce an additional state $\dagger
$. Denote by $\mathbb{T}^{\i}_{k}=\{v\in\mathbb{T}^{\i}\dvtx |v|=k\}$
the (random) collection of vertices of $\mathbb{T}^{\i}$ at distance
$k$ from the root.

The noninvasive BRW on $V\cup\{\dagger\}$ is defined together with
$(S_{v}^{\i})_{v\in\mathbb{T}^{\i}}$ on a joint probability space
such that
%
%e1.2 #&#
\begin{equation}
S_{v}^{\n}=\cases{ S_{v}^{\mathrm{aux}}, &\quad
if $S_{w}^{\i}\neq S_{v}^{\mathrm{aux}}\ \forall
w \in\mathbb{T}^{\i}_{|v|}$,
\cr
\dagger, &\quad otherwise,}
\qquad\mbox{for }v\in\mathbb{T}^{\mathrm{aux}}.
\end{equation}
We denote by $\mathbb{P}=\mathbb{P}_{\orig,x}$ the canonical
probability measure describing both processes on a same probability space.

We introduce a percolation of the family tree $\mathbb{T}^{\mathrm{aux}}$ by
declaring a vertex $v\in\mathbb{T}^{\mathrm{aux}} $ closed if and only if
$S^{\mathrm{aux}}_{v}=\dagger$. Configurations of this percolation are denoted
by $\eta_{\dagger}$. We denote by $\mathbb{T}^{\n}_{\ro}$ the
connected component of $\mathbb{T}^{\mathrm{aux}}$ containing the root $\ro$.

We say that there is \textit{coexistence} if with positive probability
both processes survive, that is, $\mathbb{P}_{\orig,x}(
|\mathbb{T}^{\n}_{\ro}|=\infty, |\mathbb{T}^{\i}|=\infty)>0$.
Using the assumptions $\mu_{\i,1}>1$ and $q_{\n}(e)>0$ together with
the strong Markov property, one sees that if $\mathbb{P}_{\orig,x}(|\mathbb{T}^{\n}_{\ro}|=\infty, |\mathbb{T}^{\i}|=\infty)>0$
holds for some $x$ then it holds for all $x\neq\orig$.

%th1.2 #&#
\begin{teo}\label{teojointsurvival1}
Let $P_{\i}$ and $P_{\n}$ transition kernels of random walks on a
infinite finitely generated group $G$ and let $\mu_{\i}$ and $\mu
_{\n}$ satisfying Assumption~\ref{ass}.
Then there is coexistence of the invasive and the noninvasive process
if $m_{\i}\rho_{\i}<1$.% and $m_{\n}\rho_{\n}\leq1$.
\end{teo}

%The proof of Theorem~\ref{teojointsurvival1} is given in Section

%re1.4 #&#
\begin{rem}[(Strongly surviving regime)]\label{remstronglysurviving}
Theorem~\ref{teojointsurvival1} states that there is always
coexistence if the invasive BRW is weakly surviving. Since we assume
the underlying random walk to be symmetric the result does not apply to
BRW on $\mathbb{Z}_{d}$. This is because Kesten's amenability
criterion implies that there is no (symmetric) weakly surviving BRW on
amenable groups (including $\mathbb{Z}^{d}$). However, on $\mathbb
{Z}^{d}$ one can show that there is no coexistence if $m_{\i}>m_{\n
}$. This can be seen by proving a shape theorem using large deviation
estimates of the underlying random walks. We refer to \cite{CP07}
where a shape theorem was established even in random environment.
However, this approach fails for groups beyond $\mathbb{Z}^{d}$ since
large deviations for random walks on groups are up to now not
sufficiently studied. Moreover, there is no reason why the shape of the
particles should be a ``convex set''; see also \cite{hueter00} and
\cite{CGM12} for results on groups with infinitely many ends. Hence,
one may ask in the flavor of Benjamini \cite{Benjamini02}: does
coexistence in the strongly surviving regime, $m_{\i}\rho_{\i}>1$
and $m_{\n}\rho_{\n}>1$, depend on the hyperbolicity of the base graph?
\end{rem}

%re1.5 #&#
\begin{rem}[(Critical case)]
We cannot treat the critical case, $m_{\i}\rho_{\i}=1$, in general
since we do not know for which groups and walks the Green function
$G(x,y|\rho_{\i}^{-1})$ decays exponentially in $d(x,y)$. However,
this is true for finite range symmetric random walks on hyperbolic
groups (see \cite{Go12}), and hence our methods also cover the
critical case in this setting.
\end{rem}

%re1.6 #&#
\begin{rem}
On groups with infinitely many ends, we have that on the event of
coexistence not every noninvasive particle has an offspring that will
be killed. This is due to the fact that the invasive BRW leaves some
neighborhoods of the boundary unvisited where the noninvasive process
may live in peace; see \cite{hueter00}~and~\cite{CGM12}. The
results in \cite{LaSe97} strongly suggests that this is still true
for Fuchsian groups and one is tempted to ask if this phenomenon holds
true for general groups. Since the shape of a single BRW is connected
to the question of coexistence for competing BRWs (see also
Remark~\ref{remstronglysurviving}), an answer to this question seems
to be related to the conjecture that the trace of a weakly surviving
BRW has infinitely many ends; see \cite{benjamini11}.
\end{rem}

%s2 #&#
\section{Preliminaries}\label{sec2}\label{secpreliminaries}

Unimodular random graphs (URGs) or stochastic homogeneous graphs have
several motivations and origins. We concentrate in this note on the
probabilistic point of view since it gives rise to the tools we are
going to use. For more details on the probabilistic viewpoints, we
refer to \cite{AL07,BC,BLS12} and to \cite{KS09}
for an introduction to the ergodic and measure theoretical origins.

One of our motivation to consider unimodular random graphs is the use
of a \emph{general} Mass-Transport Principle (MTP) which was
established in \cite{BS01} under the name of ``Intrinsic
Mass-Transport Principle'' and is basically (\ref{eqdefMTP}). It was
motivated by the fact that the Mass-Transport Principle is heavily used
in percolation theory and, therefore, lifts many results on unimodular
graphs to a more general class of graphs. In \cite{AL07}, a
probability measure on rooted graphs is called unimodular if this
general form of the MTP holds.
Another motivation to consider URGs is the fact that unimodular random
trees (URTs) can be seen as connected components in an invariant
percolation on trees; see \cite{BLS12}, Theorem 4.2, or Theorem \ref
{teounimodularinvariant} in this paper.

Let us define URGs properly. Recall that we write $(G,\orig)$ for a
graph $G=(V,E)$ with root $\orig$. A rooted isomorphism between two
rooted graphs $(G,\orig)$ and $(G',\orig')$ is an isomorphism of $G$
onto $G'$ which takes $\orig$ to $\orig'$. We denote by $\mathcal
{G}_*$ the space of isomorphism classes of rooted graphs and write
$[G,\orig]$ for the equivalence class that contains $(G,\orig)$. The
space $\mathcal{G}_{*}$ is equipped with a metric that is induced by
the~following distance between two rooted graphs $(G,\orig)$ and
$(G',\orig')$. Let $\alpha$ be the supremum of those $r>0$ such that
there exists some rooted isomorphism of~the balls of radius $\lfloor r
\rfloor$ (in graph distance) around the roots of $G$ and $G'$ and
define $d( (G,\orig), (G',\orig') )= 1/({1+\alpha})$. This metric
turns $\mathcal{G}_{*}$ into a separable and complete space. In the
same way, one defines the space $\mathcal{G}_{**}$ of isomorphism
classes of graphs with an ordered pair of distinguished vertices. A
Borel probability measure $\nu$ on $\mathcal{G}_*$ is called
unimodular if it obeys the Mass-Transport Principle: for all Borel
function $f\dvtx \mathcal{G}_{**}\to[0,\infty]$, we have
%
%e2.1 #&#
\begin{equation}
\label{eqdefMTP} \int\sum_{x\in V} f(G,\orig,x)\,d\nu
\bigl([G,\orig] \bigr)= \int\sum_{x\in V} f(G,x,\orig) \,d
\nu\bigl([G,\orig] \bigr).
\end{equation}

Observe that this definition can be extended to networks. A network is
a graph $G=(V, E)$ together with maps from $V$ and $E$ to some complete
separable metric space $\Xi$. These maps will serve as marks
(sometimes called labels) of the vertices and edges of the graph and
may a priori be unrelated. Edges are considered as directed so that
each edge is given two marks.
While the definition of the above equivalence classes for networks is
straightforward, one has to adapt the metric between two networks as
follows: $\alpha$ is chosen as the supremum of those $r>0$ such that
there is some rooted isomorphism of the balls of radius $\lfloor r
\rfloor$ around the roots of $G$ and $G'$ and such that each pair of
corresponding marks has distance at most $1/r$.
A probability measure on rooted networks is called unimodular if
equation (\ref{eqdefMTP}) holds. Realizations of these measures or
denoted as unimodular random networks (URN). Following the existing
literature, we use the notation $(G,\orig)$ as well for networks and
specify the marks of a network only when it is necessary.

Let us illustrate this definition with the very important examples of
Galton--Watson measures. Let $\mu=\{\mu_k\}_{k\in\mathbb{N}}$ be a
probability distribution on the integers. The Galton--Watson tree is
defined inductively: start with one vertex, the root of the tree. Then
the number of offspring of each particle (vertex) is distributed
according to $\mu$. Edges are between vertices and their offspring. We
denote by $\GW$ the corresponding measure on the space of rooted trees.
In this
construction, the root clearly plays a special role. For this reason,
in the unimodular Galton--Watson measure ($\UGW$) the root has a biased
distribution: the probability that the root has degree $k+1$ is
proportional to $\frac{\mu_{k}}{k+1}$. In cases where we use the $\UGW$
measure instead of the standard $\GW$ measure to define the family tree
of the BRW, we denote the BRW by UBRW.

It will be important to change the marks of a URN in a way that the
network remains unimodular. For instance, let $\xi\dvtx V\to\Xi$ be a
mark of a URN with measure~$\nu$. Let $\phi$ be a measurable map on
rooted networks that takes each network to an element of the mark space
$\Xi$. Define $\Phi$ as the map on rooted networks that takes a
network $(G,\orig)$ to another network on the same underlying graph,
but replaces the mark $\xi(x)$ by $\phi(G,x)$ for all vertices $x\in
V$. Then, by Lemma 4.1 in~\cite{BLS12}, the push forward measure
$\Phi_{*}\nu$ is also unimodular. One can also add i.i.d. marks to
existing networks. Let $\xi_{1}\dvtx V\to\Xi$ be a mark and define a new
mark $\xi(x)=(\xi_{1}(x),\xi_{2}(x))$ where the $(\xi_{2}(x))_{x\in
V}$ are realizations of i.i.d. random variables with distribution~$p$.
Denote by $\nu_{p}$ the resulting measure. Again, Lemma 4.1 in \cite
{BLS12} states that $\nu_{p}$ is unimodular.

%Suppose that $\nu$ is a unimodular measure on rooted networks. Let $
%to an element of the mark space. Define $\Phi$ to be the map on rooted
%networks that takes a network $(G,\orig)$ to another network on the
%same underlying graph, but replaces the mark at each vertex $x\in G$
%by $\phi(G,x)$. Then the push forward If instead we add an additional
%coordinate to each mark by an i.i.d.~mark, then the resulting measure
%is again unimodular.

Edge marks are maps from $V\times V\to\Xi$ and, therefore, the above
transformation can be stated analogously for edge marks. If an edge
mark $\xi$ is symmetric, that is, $\xi(x,y)=\xi(y,x)\ \forall
x,y\in V$, it can be seen as a map from $E\to\Xi$. Now, consider a
UGW-tree, add i.i.d. edge marks with distribution $q$ and denote the
resulting measure by $\UGW_{q}$. The latter marks correspond to the
steps of the BRW and we can interpret, using the definition of the UBRW
as a tree-indexed random walk, the UBRW as a URN of measure $\UGW_{q}$.

Let $\nu$ be a unimodular measure on rooted networks $(G,\orig)$ and
suppose that the mark space $\Xi$ contains a particular mark $\dagger
$. This special mark induces a natural percolation on the rooted
network: a vertex is closed if its mark equals to $\dagger$ and open
otherwise. We refer to Section~6 in \cite{AL07} for more formal
definitions and some background on percolation on URNs.

%le2.1 #&#
\begin{lem}\label{leminducedpercolation}
Let $\nu$ be a unimodular measure on rooted networks. Let $\dagger$
be a particular element of the mark space that induces a percolation.
Denote by $(C,\orig)$ the connected (marked) component containing the
origin and denote by $\nu_{\dagger}$ its corresponding measure. Then
the measure $\nu_{\dagger}$ is again a unimodular measure on rooted networks.
\end{lem}

\begin{pf}
The proof is a check of the Mass-Transport Principle (\ref
{eqdefMTP}). Define $\Phi$ as the map that takes $(G,\orig)$ to the
connected component $(C,\orig)$. The measure $\nu_{\dagger}$ is
given as the push forward measure $\Phi_{*}\nu$.
We denote by $V_{C}$ the vertex set of $(C,\orig)$ and by $V_{G}$ the
vertex set of $(G,\orig)$.
For any positive borel function $f\dvtx \mathcal{G}_{**}\to[0,\infty]$,
define its ``restrictions''
\[
f_{C}(G,x,y)=\cases{ f(C,x,y), &\quad if $x,y\in V_{C}$,
\cr
0, &\quad otherwise,}
\]
where $(C,\orig)=\Phi((G,\orig))$.
By the change of variables formula for the push forward measure, we obtain
\begin{eqnarray*}
\int \sum_{x\in V_{C}} f(C, \orig,
x) \,d\nu_{\dagger} \bigl([C,\orig] \bigr) &=& \int
\sum_{x\in\Phi([G,\orig])} f \bigl(\Phi\bigl([G,\orig] \bigr),
\orig,x \bigr) \,d\nu\bigl([G,\orig] \bigr)
\\
& = &\int
\sum_{x\in V_{G}} f_{\Phi([G,\orig])} \bigl(\Phi\bigl([G,\orig
] \bigr),\orig,x \bigr) \,d\nu\bigl([G,\orig] \bigr).
\end{eqnarray*}
Unimodularity of $\nu$ implies that the latter term is
\[
\int_{} \sum_{x\in V_{G}}
f_{\Phi([G,x])} \bigl(\Phi\bigl([G,x] \bigr),x,\orig\bigr) \,d\nu\bigl([G,
\orig] \bigr),
\]
which equals
\[
\int \sum_{x\in V_{C}} f \bigl(
\Phi\bigl([G,x] \bigr),x,\orig\bigr) \,d\nu\bigl([G,\orig] \bigr)=\int
_{} \sum_{x\in V_{C}} f(C, x, \orig) \,d
\nu_{\dagger
} \bigl([C,\orig] \bigr).
\]\upqed
\end{pf}

We make strong use of the following connection between unimodular
measures and invariant percolation.

%th2.2 #&#
\begin{teo}[(Theorem 4.2, \cite{BLS12})]\label{teounimodularinvariant}
Let $\nu$ be a probability measure on rooted networks whose underlying
graphs are trees of degree at most $d$. Then $\nu$ is unimodular iff
$\nu$ is the law of the open component of the root in a labeled
percolation on a $d$-regular tree whose law is invariant under all
automorphisms of the tree.
\end{teo}

In an invariant percolation, the probability that an edge is open is
well defined and is called the marginal of the percolation. There are
results by Adams and Lyons \cite{AL91} and H\"aggstr\"om \cite
{Hag97} that state that for invariant percolation on homogeneous trees
a sufficiently high marginal guarantees (with positive probability) the
existence of infinite clusters. We generalize this result to
``invariant percolation'' on supercritical Galton--Watson trees. This
is done by adapting the proof of Theorem~1.6 in \cite{Hag97}.

%th2.3 #&#
\begin{teo}\label{teomarginal}
Let $\UGW$ be a supercritical unimodular Galton--Watson measure of
maximal degree $d$. Then there exists some $c_{\UGW}<1$ such that for
any unimodular labeling and any particular element $\dagger$ of the
mark space the induced percolation $\UGW_{\dagger}$ assigns positive
probability to the existence of infinite clusters if the marginal is
greater than $c_{\UGW}$.
\end{teo}

\begin{pf}
Lemma~\ref{leminducedpercolation} and Theorem \ref
{teounimodularinvariant} imply that $\UGW_{\dagger}$ defines an
invariant (site)
percolation of the homogeneous tree of degree $d$.
As we have to treat two interlaced percolations, we denote by $\eta
_{\UGW}$ the configurations of the percolation induced by $\UGW$ and
by $\eta_{\dagger}$ the configurations induced by $\UGW_{\dagger}$.
From the definition of $\UGW_{\dagger}$, we have that components that
are connected in $\eta_{\dagger}$ are also connected in $\eta_{\UGW
}$. For any vertex $v$ in $\mathbb{T}_{d}$ and a given configuration
$\eta$, we write %$C_{\UGW}(e)$ for the connected component containing
%$e$ in $\eta_{\UGW}$ and
$C_{\eta}(v)$ for the connected component containing $v$ in $\eta$.
We denote by $\partial C$ the outer (vertex) boundary of a vertex set $C$.
We can now adapt the proof of Theorem 1.6 in \cite{Hag97}.
Given a configuration $\eta_{\dagger}$ we define a function $\psi
_{\eta_{\dagger}}$ on the vertex set of $\mathbb{T}_{d}$. For a
vertex $v$, denote by $v_{1},\ldots,v_{d}$ its adjacent vertices in
$\mathbb{T}^{d}$ and let
%
%e2.2 #&#
\begin{equation}
\label{eqdefpsi} \quad\psi(v)= \cases{ 1, \qquad \mbox{if }\eta_{\dagger}(v)=1 \mbox{ and }
\bigl|C_{\eta_{\dagger}}(v)\bigr|=\infty,
\cr
1, \qquad \mbox{if }\eta_{\dagger}(v)=1, \bigl|C_{\eta_{\dagger}}(v)\bigr|<\infty\mbox{ and } \displaystyle\frac{ | C_{\eta_{\dagger
}}(v)|}{|\partial C_{ \eta
_{\dagger}}(v)| }\geq K,
\vspace*{5pt}\cr
0, \qquad \mbox{if }\eta_{\dagger}(v)=1, \bigl|C_{\eta_{\dagger}}(v)\bigr|<\infty\mbox{ and }
\displaystyle\frac{ | C_{\eta_{\dagger}}(v)|}{|\partial C_{\eta_{\dagger}}(v)| }<K,
\cr
\displaystyle 1 + \sum_{i=1}^{d}
f(v_{i}),
\cr
\hspace*{32.5pt}\mbox{if }\eta_{\dagger}(v)=0,}
\end{equation}
where
\[
f(w) = \cases{ \displaystyle\frac{ | C_{\eta_{\dagger}}(w)|}{|\partial
C_{\eta_{\dagger}}(w)|}, &\quad if $\displaystyle | C_{\eta_{\dagger}}(w)|<
\infty\mbox{ and }\frac{
| C_{\eta_{\dagger}}(w)|}{|\partial C_{\eta_{\dagger}}(w)| }<K$,
\vspace*{5pt}\cr
0, &\quad otherwise,}
\]
for some positive constant $K$ to be chosen later.
We can now, as in the proof Theorem 1.6 in \cite{Hag97}, interpret
$\psi$ as a distribution of mass over the vertices. Originally, every
vertex has mass $1$. For vertices $v$ in an infinite cluster or
vertices $v$ in finite clusters with $\frac{ | C_{\eta_{\dagger
}}(v)|}{|\partial C_{\eta_{\dagger}}(v)| }\geq K$, the mass in $v$
remains unchanged. If $v$ is in a finite cluster such that $\frac{ |
C_{\eta_{\dagger}}(v)|}{|\partial C_{\eta_{\dagger}}(v)| }< K$,
then $v$ distributes its mass equally to the closed vertices incident
to $C_{\eta_{\dagger}}(v)$.
If $\eta_{\dagger}(v)=0$, then $v$ receives additional mass from the
distributing vertices. For two vertices $v$ and $w$, we write
$\triangle\psi(v,w)$
for the flow of mass from $v$ to $w$ (using the above interpretation)
and obtain
\[
\psi(v)=1+\sum_{w} \triangle\psi(w,v).
\]
Consider random configurations $X_{\UGW}$ and $X_{\dagger}$ that are
distributed according to $\UGW$ and $\UGW_{\dagger}$. Since $\UGW
_{\dagger}$ is unimodular, we have for any pair of vertices $v$ and
$w$ that
\[
\mathbb{E}_{\UGW_{\dagger}} \bigl[\triangle\psi(v,w) \bigr]=0.
\]
Since $\psi$ is bounded, we obtain that
%
%e2.3 #&#
\begin{equation}
\label{eqmtp} \mathbb{E}_{\UGW_{\dagger}} \bigl[\psi(v) \bigr] = 1 +
\sum
_{w} \mathbb{E}_{\UGW_{\dagger}} \bigl[\triangle
\psi(w,v) \bigr]=1.
\end{equation}
For\vspace*{-3pt} the sake of typesetting, we write $\{ \frac{|C_{X_\dagger
}(v)|}{|\partial C_{X_\dagger}(v)| }\geq K\}$ for the event $\{
|C_{X_\dagger}(v)|<\infty, \frac{|C_{X_\dagger}(v)|}{|\partial
C_{X_\dagger}(v)| }\geq K \}$.
Using equation (\ref{eqmtp})\vspace*{1pt} with the definition of $\psi$ in
equation (\ref{eqdefpsi}), we obtain that
\[
\UGW_{\dagger} \bigl( X_{\dagger}(v)=1, \bigl|C_{X_\dagger}(v)\bigr|=\infty
\bigr)
\]
is greater or equal than
\[
1-\UGW_{\dagger} \biggl(X_{\dagger}(v)=1, \frac{|C_{X_\dagger
}(v)|}{|\partial C_{X_\dagger}(v)| }\geq K
\biggr) - c_{K} \UGW_{\dagger} \bigl(X_{\dagger}(v)=0
\bigr),
\]
where $c_{K}=1+dK$.
In order to adjust the value of $K$ recall the following.
The anchored (vertex) isoperimetric constant for a graph $G$ is defined as
\[
\i(G,v)=\inf_{S\ni v} \frac{|\partial S|}{|S|},
\]
where $S$ ranges over all connected vertex sets containing a fixed
vertex $v$. Note that $\i(G,v)$ does not depend on the choice of the
edge $v$. Corollory 1.3 in \cite{CP04} states %, however formulated
%for the vertex isoperimetric constant, that for a supercritical
%Galton-Watson tree $\T$ we have
that $\i(\mathbb{T},v)>0$ a.s. on the event that $\mathbb{T}$ is
infinite. Now, since
\[
\UGW_{\dagger} \biggl(X_{\dagger}(v)=1,\frac{|C_{X_{\dagger
}}(v)|}{|\partial C_{X_{\dagger}}(v)| }\geq K
\biggr)
\]
is bounded above by
\[
\UGW\bigl(\i(\mathbb{T})^{-1}>K, \bigl|C_{X_\UGW}(v)\bigr|=\infty\bigr),
\]
we can choose $K$ sufficiently large such that
\[
\UGW_{\dagger} \biggl(X_{\dagger}(v)=1, \frac{|C_{X_{\dagger
}}(v)|}{|\partial C_{X_{\dagger}}(v)| }\geq K
\biggr)< \UGW\bigl(\bigl|C_{X_\UGW}(v)\bigr|=\infty\bigr).
\]
Eventually, there exists some constant $c>0$ such that
\[
\UGW_{\dagger} \bigl( \bigl|C_\dagger(v)\bigr|=\infty| \bigl|C_{X_\UGW}(v)\bigr|=
\infty\bigr) > c- \frac{c_{K} \UGW_{\dagger}(X_{\dagger}(v)=0))}{ \UGW
(|C_{X_\UGW}(v)|=\infty)}.
\]
Hence, choosing the marginal $ \UGW_{\dagger}( X_{\dagger}(v)=1)$
sufficiently high assures that $\UGW_{\dagger}( |C_{X_\dagger
}(v)|=\infty)>0$.
\end{pf}

%
%re2.1 #&#
\begin{rem}
An inspection of the proof above reveals that Theorem \ref
{teomarginal} holds true for unimodular measures on rooted networks
whose underlying graphs are trees of bounded degree and that give
positive weight to infinite networks such that almost all infinite
realizations have positive anchored isoperimetric constant.
\end{rem}

%s3 #&#
\section{Truncated BRW---proof of Theorem \texorpdfstring{\protect\ref{teotruncated}}{1.1}}\label{sec3}
Since the case $m>1/\rho(P)$ was proven in \cite{Zucca11} let us
assume in the following that $m\leq1/\rho(P)$.
For the case $\mu_{0}=0$, the proof is essentially given in the proof
of Theorem 3.1 in \cite{benjamini11}. We give a concise proof using
the results of \cite{BLS12} and Theorem~\ref{teomarginal}.
Moreover, we hope that the example of truncated BRW will serve as an
introduction of our approach ``interacting growth process and invariant
percolation,'' and hence is useful for a better understanding of the
proof of Theorem~\ref{teojointsurvival1}. Our approach consists of
two steps: an adaptation of the model such that the family tree is a
URT and the control of the marginal of the corresponding invariant percolation.

%s3.1 #&#
\subsection{Adapting the model}\label{sec3.1}
The aim is to identify an invariant percolation (or unimodular measure)
of the family tree. Since the percolation induced by $\dagger$ is in
not an invariant percolation, there is need for a reformulation of our
problem. We will define a new process in a way such that vertices that
were visited more than $N$ times become ``deadly'' for all instances of
times. In other words, if $x$ is a vertex of the base graph such that
$|\{v\dvtx S_{v}=x\}|> N$, then we set $S_{v}^{\mathrm{new}}=\dagger$ for all $v$
such that $S_{v}=x$. More formally, let $(\mathbb{T},\ro)$ be the
labeled UGW-tree (the BRW) and define
%
%e3.1 #&#
\begin{equation}
\phi(\mathbb{T},x)=\cases{ \dagger, &\quad$\bigl|\{v\dvtx S_{v}=x\}\bigr|> N$,
\vspace*{3pt}\cr
\bullet, &\quad$\bigl|\{v\dvtx S_{v}=x\}\bigr|\leq N$.}
\end{equation}
The corresponding push forward measure $\Phi_{*} \UGW_{q}$ is again
unimodular; see Lemma 4.1 in \cite{BLS12}.

%s3.2 #&#
\subsection{Control of the marginal}\label{sec3.2}
The underlying BRW is supposed to be weakly surviving, that is,
$\mathbb{P}( |\{v\dvtx S_{v}=S_{\ro}\}|<\infty)=1$. Hence, we can apply
Theorem~\ref{teomarginal} and choose $N_{u}$ sufficiently large such
that the marginal $\mathbb{P}( |\{v\dvtx S_{v}=S_{\ro}\}|\leq N_{u})$ is
sufficiently high. This guarantees that with positive probability the
cluster containing $\ro$ is infinite and that the process
$S_{v}^{\mathrm{new}}$ survives with positive probability. Since $S_{v}^{\mathrm{new}}$
is stochastically dominated by the truncated BRW, we obtain that
$\BRW_{N}$ survives with positive probability for sufficiently large
$N$. This yields, together with the monotonicity of the model, the
existence of a critical value $N_{c}$ given in Theorem~\ref{teotruncated}.

%s4 #&#
\section{Competing BRWs---proof of Theorem \texorpdfstring{\protect\ref{teojointsurvival1}}{1.2}}\label{sec4}

We proceed in two steps as in the previous section. In this section, we
suppose that $\rho_{\i} m_{\i}<1$ and assume without loss of
generality that $\rho_{\n}m_{\n}< 1$.
%s4.1 #&#
\subsection{Adapting the model}\label{sec4.1}\label{subseccompeting1}
The family tree of the noninvasive process is in general not a URT.
We invite the reader to the following informal description of the situation.

Let us start both processes in neighboring sites, then the offspring of
the starting particles are very likely to be killed by those of the
invasive process. However, if we consider some noninvasive particle
very late in the genealogical process, then given the fact that the
particle exists (or is alive), one might expect that its ancestors
never have been very close to invasive particles. Hence, the chances of
its children to survive are high as well. As a conclusion, we have to
adapt the invasive process in a way that every particle of the
auxiliary process (of the noninvasive process) has the same probability
to encounter an invasive particle. For this purpose, we will not start
just one invasive process but infinitely many.

In the following, we describe a first approach that gives the right idea
but does not lead to a good control of the marginal.
First of all, there is a natural mapping $v\mapsto S_{v}$ from the
family tree of the auxiliary process to the base graph that we denote
by $\Psi$. Now, on the base graph we start infinitely many independent
BRWs according to $(\mu_{\i}, P_{\i})$ as follows. Let $N\in\mathbb
{N}$ (to be chosen later) and start independent copies of invasive BRWs
on each $x$ with $|\Psi^{-1}(x)|=N$. Here, it is important that the
underlying BRW survives weakly; otherwise the latter set would be
empty. Using these BRWs, we define a random labeling of the base graph
$G$: a vertex is labeled $\dagger$ if it is visited by some invasive
particle at some time and $\bullet$ otherwise. In \cite
{benjamini11}, it was shown that the trace of a (weakly surviving) BRW
is a URG, and moreover that the above labeling defines a URN. We use
now the map $\Psi$ to retrieve this labeling; label a vertex $v\in
\mathbb{T}^{\mathrm{aux}}$ with $\dagger$ if $\Psi(v)$ is labeled by $\dagger
$ and label it with $\bullet$ otherwise. Each of the above steps is
invariant under rerooting and so is the new labeled version of $\mathbb
{T}^{\mathrm{aux}}$. Finally, due to Lemma~\ref{leminducedpercolation}, the
connected component of $\mathbb{T}^{\mathrm{aux}}$ with respect to the
percolation induced by $\dagger$ is a URT. It remains to prove that
the noninvasive BRW survives with positive probability when being
confronted with an infinity of invasive BRWs. This would imply
coexistence of the two original processes, since coexistence does not
depend on the starting position of the processes.

%s4.2 #&#
\subsection{Control of the marginal}\label{sec4.2}
In general, it is not possible to control the marginal of the invariant
percolation above. In fact, we need a better control of the ``number''
of invasive processes. Denote by $\mathcal{B}(n,\orig)=\{x\dvtx d(\orig,x)\leq n\}$ the ball of radius $n$ around the origin $\orig$ and
denote by $\mathcal{S}(\orig,n)=\{x\dvtx d(\orig,x)= n\}$ the
corresponding sphere.
The growth rate $g$ of the group $\Gamma$ is defined as $g=\lim_{n\to
\infty} \frac{1}n \log(|\mathcal{B}(n,\orig)|)$.

As the underlying random walks are supposed to be symmetric random
walks we have (see \cite{woess}, Lemma 8.1) that $p^{(n)}(x,y)\leq
\rho^{n}$ for all $x,y\in\Gamma$ and all $n\in\mathbb{N}$. Two
consequences of this fact on BRWs are given in the following lemma.

%le4.1 #&#
\begin{lem}\label{lemcontrol}
Let $(\mu,P)$ be a BRW on a nonamenable Cayley graph with $\rho m<1$. Then:
\begin{longlist}[(2)]
\item[(1)] $G(x,y | m):=\sum_{n=0}^{\infty} p^{(n)}(x,y) m^{n}\leq
(m\rho)^{d(x,y)}/{(1-\rho m)}$;
\item[(2)] there exists some constant $\ell$ such that
\[
\limsup_{n\to\infty} \mathbb{E} \bigl[ \bigl| \bigl\{v\dvtx S_{v}
\in\mathcal{B}(\orig,n) \bigr\}\bigr| \bigr]/ m^{\ell n}=0.
\]
\end{longlist}
\end{lem}

\begin{pf}
(1) Since the random walk is nearest neighbor, that is,
$\operatorname{supp}(q)=S$, we have that
\[
G(x,y|m)=\sum_{n=d(x,y)}^{\infty}
p^{(n)}(x,y) m^{n} \leq\sum_{n=0}^{\infty}
(\rho m)^{n+d(x,y)} \leq(m\rho)^{d(x,y)} \frac{1}{1-\rho m}.
\]

(2)
%We have\begin{eqnarray*} \E[|\{v\dvtx S_{v}\in\B(\orig,n) \}| &=&
%&\leq& \sum_{k\geq0} \sum_{y\in\B(\orig,n)} m^{k} p^{(k)}(\orig,y)\cr
%&\leq& C \sum_{k\geq0} m^{k} \rho^{k} g^{n} \leq C g^{n}
Denote $R_{n}=\inf\{k\geq0\dvtx S_{v}\notin\mathcal{B}(\orig,n)\ \forall
|v|\geq k\}$. In the following, denote by $C$ a constant
that is always chosen sufficiently large and may change from formula to formula.
For some $b>0$ (to be chosen in a moment), we obtain using the Markov inequality
\begin{eqnarray*}
\mathbb{P}(R_{n}> b n) &=& \mathbb{P} \bigl(\exists v\dvtx |v|\geq b
n\dvtx S_{v}\in\mathcal{B}(\orig, n) \bigr)
\\
&\leq& \sum
_{k\geq b n}\sum_{y\in\mathcal{B}(\orig, n)} m^{k}
p^{(k)}(\orig, y)
\\
&\leq& C \sum_{k\geq b n}
m^{k} \rho^{k} g^{n} \leq C \bigl(g (m
\rho)^{b} \bigr)^{n}.
\end{eqnarray*}
Hence, we can choose $b$ sufficiently large such that the latter
probability is summable and $\limsup R_{n}/n \leq b$ by the lemma of
Borel--Cantelli. Finally, for $n$ sufficiently large
\begin{eqnarray*}
\mathbb{E} \bigl[\bigl| \bigl\{v\dvtx S_{v}\in\mathcal{B}(\orig,n) \bigr\}\bigr|
\bigr] &\leq& \mathbb{E} \bigl[\bigl|\bigl\{v\dvtx|v|\leq R_{n}\bigr\}\bigr| \bigr]
\\
& \leq&
\mathbb{E} \bigl[ \bigl| \bigl\{v\dvtx|v|\leq(b+1)n \bigr\}\bigr| \bigr]
\\
& \leq&
\frac{m^{(b+1)n+1}-1}{m-1} \leq m^{(b+1)n+1}, %\leq\E[\{v\dvtx|v|\leq
\end{eqnarray*}
which yields the result for some sufficiently large $\ell$.
\end{pf}

The first part of Lemma~\ref{lemcontrol} is used to control each of
the invasive processes and the second part to adjust the ``number'' of
these invasive processes.
In order to start with the adjustment, let us consider a noninvasive
process with less branching. For any constant $\gamma\in(0,1]$, to be
chosen later, we define the truncated Galton--Watson process by
\[
\mu^{(\gamma)}_{k}= \cases{ \displaystyle\gamma\mu^{}_{\n,k},
&\quad for $k\geq2$,
\cr
\displaystyle\mu_{\n,1} + (1-\gamma) \sum
_{k=2}^{\infty} \mu^{}_{\n,k}, &
\quad for $k= 1$,
\cr
\mu_{\n,0},&\quad for $k=0$}
\]
and denote its mean by $m_{\gamma}$. This construction is made to
ensure two main properties: $m_{\gamma}\searrow1-\mu_{\n,0}\leq1$
as $\gamma\searrow0$ and $\mu^{(\gamma_{1})}_{k}<\mu^{(\gamma
_{2})}_k$ for all $\gamma_{1}<\gamma_{2}$ and $k\geq2$. This latter
property allows to construct a natural coupling of the original and the
``$\gamma$-processes.'' Hence, denote by $S^{\gamma}$ the BRW
corresponding to the family tree $\mathbb{T}^{\gamma}$. Due to the
coupling, it remains to show that the ``$\gamma$-process'' has
positive probability of survival for some $\gamma>0$. Let $\gamma
_{c}$ be such that $m_{\gamma_{c}}=1$.

Recall the definition of $\Psi$ and $\dagger$ in Section~\ref
{subseccompeting1} and start independent copies of invasive BRWs on
each $x$ with $|\Psi^{-1}(x)|=N$. (The constant $N$ is still to be chosen.)
We will also denote by $\mathbb{P}$ the probability measure describing
the noninvasive process together with the infinite number of invasive processes.
Denote by $A\subset G$ the (random) set where invasive processes are
started. Since $A\cap\mathcal{S}(\orig,n)\subset\{x\in\mathcal
{S}(\orig,n): \exists v\in S_{v}^{\gamma}=x\}$ and vertices in $A$
are labeled by $\dagger$, we have
\[
\label{eqstar} \sum_{x\in\mathcal{S}(\orig,n)} \mathbb{P}(x\in A)\leq
\mathbb{E} \bigl[\bigl| \bigl\{v\dvtx S_{v}^{\gamma}\in\mathcal{S}(
\orig,n), \xi(v)=\dagger\bigr\}\bigr| \bigr],
\]
where $\xi(v)$ denotes the mark of the vertex $v$ induced by $\Psi$.
For $x\in A$, we denote by $S^{\i,x}_{v}$ the invasive BRW started in
$x$ with family tree $\mathbb{T}^{\i,x}$.
Due to Lemma~\ref{lemcontrol} for any $\gamma\in(\gamma_{c},1)$,
there exists some constants $C_{\gamma}$ and $\ell_{\gamma}$ such
that $\mathbb{E}[|\{v\dvtx S_{v}^{\gamma}\in\mathcal{S}(\orig,n)\}
|]\leq C_{\gamma} m_{\gamma}^{\ell_{\gamma} n}$.
Since the trace of the BRW is unimodular, we have that there exists a
constant $C_{N}\to0$ (as $N\to\infty$) such that
\[
\mathbb{E} \bigl[\bigl| \bigl\{v\dvtx S_{v}^{\gamma}\in\mathcal{S}(
\orig,n), \xi(v)=\dagger\bigr\}\bigr| \bigr] \leq C_{N} C_{\gamma}
m_{\gamma}^{\ell_{\gamma} n}.
\]
Moreover, the proof of Lemma~\ref{lemcontrol} gives that the constant
$\ell_{\gamma}$ can be chosen uniform with respect to $\gamma$ since
there is a natural coupling for the last exit times $R_{n}$ of
different ``$\gamma$-processes.'' Hence, there exists some constant
$\ell$ such that for all $\gamma\in(\gamma_{c},1]$
\[
\label{eq2star} \sum_{x\in\mathcal{S}(\orig,n)} \mathbb{P}(x\in A) \leq
C_{N} C_{\gamma} m_{\gamma}^{\ell n}.
\]
Using this together with a union bound and part (1) of Lemma \ref
{lemcontrol}, we obtain
\begin{eqnarray*}
\mathbb{P} \bigl( \xi(\ro)=\dagger\bigr) &\leq& \mathbb{P} \bigl
(\exists x\in
A, \exists v\in\mathbb{T}^{\i,x}\dvtx S^{\i,x}_{v}=S_{r}^{\gamma
}
\bigr)
\\
& \leq& \sum_{x\in G} \mathbb{P} \bigl(x\in
A, \exists v\in\mathbb{T}^{\i,x}\dvtx S^{\i,x}_{v}=S_{r}^{\gamma}
\bigr)
\\
& =& \sum_{n=0}^{\infty} \sum
_{x\in\mathcal{S}(\orig, n)} \mathbb{E} \bigl[\bigl| \bigl\{
v\dvtx S^{\i,x}_{v}=S^{\gamma}_{\ro} \bigr\}\bigr| x\in A
\bigr] \mathbb{P}(x\in A)
\\
& \leq& \sum_{n=0}^{\infty}
\frac{1}{1-\rho_{\i}m_{\i}} (m_{\i
} \rho_{\i})^{n} \sum
_{x\in\mathcal{S}(\orig, n)} \mathbb{P}(x\in A)
\\
& \leq& \sum
_{n=0}^{\infty} \frac{1}{1-\rho_{\i}m_{\i}} (m_{\i
}
\rho_{\i})^{n} C_{N} C_{\gamma}
m_{\gamma}^{\ell n}.%\cr
% & \leq& \sum_{n=0}^{\infty} \frac1{m_{\i}\rho_{\i}}C_{N} C_{\gamma}
%(m_{\gamma}^{\ell}) (m_{\i}\rho_{\i})^{n}
\end{eqnarray*}
We can choose $\gamma\in(\gamma_{c},1]$ sufficiently small such that
$m_{\gamma}^{\ell} m_{\i}\rho_{\i}<1$. Let $c_{\UGW_{\gamma}}$
be the constant from Theorem~\ref{teomarginal} for the Galton--Watson
with offspring distribution~$\mu^{(\gamma)}$. Now, choose $N$
sufficiently large (which makes $C_{N}$ sufficiently small) such that
the marginal $\mathbb{P}( \Xi(\ro)\neq\dagger)>c_{\UGW_{\gamma
}}$. This in turn implies that the noninvasive BRW with offspring
distribution $\mu^{(\gamma)}$ survives with positive probability if
confronted with an infinite number of invasive BRWs. Hence, for some
$\gamma_{c'}\in(\gamma_{c},1]$ there is coexistence of one invasive
and one noninvasive BRW since coexistence does not depend on the choice
of the starting positions of the processes. Eventually, using the
monotonicity in $\gamma$ a standard coupling argument implies
coexistence for all $\gamma\in[\gamma_{c'},1]$.

\section*{Acknowledgments}
The idea behind this paper came up during the focussed meeting
``Branching Random Walks and Related Topics'' in Graz. The author
thanks TU-Graz and RGLIS for organization and funding. Moreover, the
author appreciated the hospitality of the Department of Mathematical
Structure Theory of TU-Graz during the write up.
The author thanks Itai Benjamini for helpful comments on a first
version of this note and is grateful to the referee for various
valuable comments.


%suskaldyti doi

% imsref loaded by linak, 2014-03-20 12:13:52
%

\printaddresses


\begin{thebibliography}{28}
% pybtex-1.01. Style name=ims, version=2.7, label_style=nolabel,
%sorting_style=complex, cfg=None, language=None.
%b1 ###
%b1 #&#
\bibitem{AL91}
%
\begin{barticle}[mr]
\bauthor{\bsnm{Adams},~\bfnm{Scot}\binits{S.}} \AND
\bauthor{\bsnm{Lyons},~\bfnm{Russell}\binits{R.}}
(\byear{1991}).
\btitle{Amenability, {K}azhdan's property and percolation for trees,
groups and equivalence relations}.
\bjournal{Israel J. Math.}
\bvolume{75}
\bpages{341--370}.
\bid{doi={10.1007/BF02776032}, issn={0021-2172}, mr={1164598}}
\end{barticle}
%
\bptok{imsref}%
% NOT OUTPUTED:
% issn = 0021-2172
% url = http://dx.doi.org/10.1007/BF02776032
% number = 2-3
% coden = ISJMAP
% fjournal = Israel Journal of Mathematics
\endbibitem

%b2 ###
%b2 #&#
\bibitem{AL07}
%
\begin{barticle}[mr]
\bauthor{\bsnm{Aldous},~\bfnm{David}\binits{D.}} \AND
\bauthor{\bsnm{Lyons},~\bfnm{Russell}\binits{R.}}
(\byear{2007}).
\btitle{Processes on unimodular random networks}.
\bjournal{Electron. J. Probab.}
\bvolume{12}
\bpages{1454--1508}.
\bid{doi={10.1214/EJP.v12-463}, issn={1083-6489}, mr={2354165}}
\end{barticle}
%
\bptok{imsref}%
% NOT OUTPUTED:
% issn = 1083-6489
% url = http://dx.doi.org/10.1214/EJP.v12-463
% fjournal = Electronic Journal of Probability
\endbibitem

%b3 ###
%b3 #&#
\bibitem{Amarasekare}
%
\begin{barticle}[author]
\bauthor{\bsnm{Amarasekare},~\bfnm{P.}\binits{P.}}
(\byear{2003}).
\btitle{Competitive coexistence in spatially structured environments:
A~synthesis}.
\bjournal{Ecology Letters}
\bvolume{6}
\bpages{1109--1122}.
\end{barticle}
%
\bptok{imsref}%
\endbibitem

%b4 ###
%b4 #&#
\bibitem{Benjamini02}
%
\begin{barticle}[mr]
\bauthor{\bsnm{Benjamini},~\bfnm{Itai}\binits{I.}}
(\byear{2002}).
\btitle{Survival of the weak in hyperbolic spaces, a remark on
competition and geometry}.
\bjournal{Proc. Amer. Math. Soc.}
\bvolume{130}
\bpages{723--726 (electronic)}.
\bid{doi={10.1090/S0002-9939-01-06077-4}, issn={0002-9939}, mr={1866026}}
\end{barticle}
%
\bptok{imsref}%
% NOT OUTPUTED:
% issn = 0002-9939
% url = http://dx.doi.org/10.1090/S0002-9939-01-06077-4
% number = 3
% coden = PAMYAR
% fjournal = Proceedings of the American Mathematical Society
\endbibitem

%b5 ###
%b5 #&#
\bibitem{BC}
%
\begin{barticle}[mr]
\bauthor{\bsnm{Benjamini},~\bfnm{Itai}\binits{I.}} \AND
\bauthor{\bsnm{Curien},~\bfnm{Nicolas}\binits{N.}}
(\byear{2012}).
\btitle{Ergodic theory on stationary random graphs}.
\bjournal{Electron. J. Probab.}
\bvolume{17}
\bpages{1--20}.
\bid{doi={10.1214/EJP.v17-2401}, issn={1083-6489}, mr={2994841}}
\end{barticle}
%
\bptok{imsref}%
% NOT OUTPUTED:
% issn = 1083-6489
% url = http://dx.doi.org/10.1214/EJP.v17-2401
% fjournal = Electronic Journal of Probability
\endbibitem

%b6 ###
%b6 #&#
\bibitem{BLS12}
%
\begin{bmisc}[author]
\bauthor{\bsnm{Benjamini},~\bfnm{I.}\binits{I.}},
\bauthor{\bsnm{Lyons},~\bfnm{R.}\binits{R.}} \AND
\bauthor{\bsnm{Schramm},~\bfnm{O.}\binits{O.}}
(\byear{2013}).
\bhowpublished{Unimodular random trees.
\textit{Ergodic Theory Dynam. Systems}. DOI:\doiurl{10.1017/etds.2013.56}.}
\end{bmisc}
%
\bptok{imsref}%
\endbibitem

%b7 ###
%b7 #&#
\bibitem{benjamini11}
%
\begin{barticle}[mr]
\bauthor{\bsnm{Benjamini},~\bfnm{Itai}\binits{I.}} \AND
\bauthor{\bsnm{M{\"u}ller},~\bfnm{Sebastian}\binits{S.}}
(\byear{2012}).
\btitle{On the trace of branching random walks}.
\bjournal{Groups Geom. Dyn.}
\bvolume{6}
\bpages{231--247}.
\bid{doi={10.4171/GGD/156}, issn={1661-7207}, mr={2914859}}
\end{barticle}
%
\bptok{imsref}%
% NOT OUTPUTED:
% issn = 1661-7207
% url = http://dx.doi.org/10.4171/GGD/156
% number = 2
% fjournal = Groups, Geometry, and Dynamics
\endbibitem

%b8 ###
%b8 #&#
\bibitem{benjamini94}
%
\begin{barticle}[mr]
\bauthor{\bsnm{Benjamini},~\bfnm{Itai}\binits{I.}} \AND
\bauthor{\bsnm{Peres},~\bfnm{Yuval}\binits{Y.}}
(\byear{1994}).
\btitle{Markov chains indexed by trees}.
\bjournal{Ann. Probab.}
\bvolume{22}
\bpages{219--243}.
\bid{issn={0091-1798}, mr={1258875}}
\end{barticle}
%
\bptok{imsref}%
% NOT OUTPUTED:
% issn = 0091-1798
% url =
%http://links.jstor.org/sici?sici=0091-1798(199401)22:1<219:MCIBT>2.0.CO;2-U&origin=MSN
% number = 1
% coden = APBYAE
% fjournal = The Annals of Probability
\endbibitem

%b9 ###
%b9 #&#
\bibitem{BS01}
%
\begin{barticle}[mr]
\bauthor{\bsnm{Benjamini},~\bfnm{Itai}\binits{I.}} \AND
\bauthor{\bsnm{Schramm},~\bfnm{Oded}\binits{O.}}
(\byear{2001}).
\btitle{Recurrence of distributional limits of finite planar graphs}.
\bjournal{Electron. J. Probab.}
\bvolume{6}
\bpages{no. 23, 13 pp. (electronic)}.
\bid{doi={10.1214/EJP.v6-96}, issn={1083-6489}, mr={1873300}}
\end{barticle}
%
\bptok{imsref}%
% NOT OUTPUTED:
% issn = 1083-6489
% url = http://dx.doi.org/10.1214/EJP.v6-96
% fjournal = Electronic Journal of Probability
\endbibitem

%b10 ###
%b10 #&#
\bibitem{BPZ07}
%
\begin{barticle}[mr]
\bauthor{\bsnm{Bertacchi},~\bfnm{Daniela}\binits{D.}},
\bauthor{\bsnm{Posta},~\bfnm{Gustavo}\binits{G.}} \AND
\bauthor{\bsnm{Zucca},~\bfnm{Fabio}\binits{F.}}
(\byear{2007}).
\btitle{Ecological equilibrium for restrained branching random walks}.
\bjournal{Ann. Appl. Probab.}
\bvolume{17}
\bpages{1117--1137}.
\bid{doi={10.1214/105051607000000203}, issn={1050-5164}, mr={2344301}}
\end{barticle}
%
\bptok{imsref}%
% NOT OUTPUTED:
% issn = 1050-5164
% url = http://dx.doi.org/10.1214/105051607000000203
% number = 4
% fjournal = The Annals of Applied Probability
\endbibitem

%b11 ###
%b11 #&#
\bibitem{BZ09}
%
\begin{barticle}[mr]
\bauthor{\bsnm{Bertacchi},~\bfnm{Daniela}\binits{D.}} \AND
\bauthor{\bsnm{Zucca},~\bfnm{Fabio}\binits{F.}}
(\byear{2009}).
\btitle{Approximating critical parameters of branching random walks}.
\bjournal{J. Appl. Probab.}
\bvolume{46}
\bpages{463--478}.
\bid{doi={10.1239/jap/1245676100}, issn={0021-9002}, mr={2535826}}
\end{barticle}
%
\bptok{imsref}%
% NOT OUTPUTED:
% issn = 0021-9002
% url = http://dx.doi.org/10.1239/jap/1245676100
% number = 2
% coden = JPRBAM
% fjournal = Journal of Applied Probability
\endbibitem

%b12 ###
%b12 #&#
\bibitem{Blair-Stahn10}
%
\begin{bmisc}[author]
\bauthor{\bsnm{Blair-Stahn},~\bfnm{N.~D.}\binits{N.~D.}}
(\byear{2010}).
\bhowpublished{First passage percolation and competition models.
Preprint. Available at \arxivurl{arXiv:1005.0649}.}
\end{bmisc}
%
\bptok{imsref}%
\endbibitem

%b13 ###
%b13 #&#
\bibitem{CGM12}
%
\begin{barticle}[mr]
\bauthor{\bsnm{Candellero},~\bfnm{Elisabetta}\binits{E.}},
\bauthor{\bsnm{Gilch},~\bfnm{Lorenz~A.}\binits{L.~A.}} \AND
\bauthor{\bsnm{M{\"u}ller},~\bfnm{Sebastian}\binits{S.}}
(\byear{2012}).
\btitle{Branching random walks on free products of groups}.
\bjournal{Proc. Lond. Math. Soc. (3)}
\bvolume{104}
\bpages{1085--1120}.
\bid{doi={10.1112/plms/pdr060}, issn={0024-6115}, mr={2946082}}
\end{barticle}
%
\bptok{imsref}%
% NOT OUTPUTED:
% issn = 0024-6115
% url = http://dx.doi.org/10.1112/plms/pdr060
% number = 6
% fjournal = Proceedings of the London Mathematical Society. Third
%Series
\endbibitem

%b14 ###
%b14 #&#
\bibitem{CP04}
%
\begin{barticle}[mr]
\bauthor{\bsnm{Chen},~\bfnm{Dayue}\binits{D.}} \AND
\bauthor{\bsnm{Peres},~\bfnm{Yuval}\binits{Y.}}
(\byear{2004}).
\btitle{Anchored expansion, percolation and speed}.
\bjournal{Ann. Probab.}
\bvolume{32}
\bpages{2978--2995}.
\bid{doi={10.1214/009117904000000586}, issn={0091-1798}, mr={2094436}}
\end{barticle}
%
\bptok{imsref}%
% NOT OUTPUTED:
% issn = 0091-1798
% url = http://dx.doi.org/10.1214/009117904000000586
% number = 4
% coden = APBYAE
% fjournal = The Annals of Probability
\endbibitem

%b15 ###
%b15 #&#
\bibitem{CP07}
%
\begin{barticle}[mr]
\bauthor{\bsnm{Comets},~\bfnm{Francis}\binits{F.}} \AND
\bauthor{\bsnm{Popov},~\bfnm{Serguei}\binits{S.}}
(\byear{2007}).
\btitle{Shape and local growth for multidimensional branching random
walks in random environment}.
\bjournal{ALEA Lat. Am. J. Probab. Math. Stat.}
\bvolume{3}
\bpages{273--299}.
\bid{issn={1980-0436}, mr={2365644}}
\end{barticle}
%
\bptok{imsref}%
% NOT OUTPUTED:
% issn = 1980-0436
% fjournal = ALEA. Latin American Journal of Probability and
%Mathematical Statistics
\endbibitem

%b16 ###
%b16 #&#
\bibitem{gantert04}
%
\begin{barticle}[mr]
\bauthor{\bsnm{Gantert},~\bfnm{N.}\binits{N.}} \AND
\bauthor{\bsnm{M{\"u}ller},~\bfnm{S.}\binits{S.}}
(\byear{2006}).
\btitle{The critical branching {M}arkov chain is transient}.
\bjournal{Markov Process. Related Fields}
\bvolume{12}
\bpages{805--814}.
\bid{issn={1024-2953}, mr={2284404}}
\bptnote{check year}%
\end{barticle}
%
\bptok{imsref}%
% NOT OUTPUTED:
% issn = 1024-2953
% number = 4
% fjournal = Markov Processes and Related Fields
\endbibitem

%b17 ###
%b17 #&#
\bibitem{Go12}
%
\begin{bmisc}[author]
\bauthor{\bsnm{Gouezel},~\bfnm{S.}\binits{S.}}
(\byear{2014}).
\bhowpublished{Local limit theorem for symmetric random walks in
Gromov-hyperbolic groups.
\textit{J. Amer. Math. Soc.} DOI:\doiurl{10.1090/S0894-0347-2014-00788-8}.}
\end{bmisc}
%
\bptok{imsref}%
\endbibitem

%b18 ###
%b18 #&#
\bibitem{Hag97}
%
\begin{barticle}[mr]
\bauthor{\bsnm{H{\"a}ggstr{\"o}m},~\bfnm{Olle}\binits{O.}}
(\byear{1997}).
\btitle{Infinite clusters in dependent automorphism invariant
percolation on trees}.
\bjournal{Ann. Probab.}
\bvolume{25}
\bpages{1423--1436}.
\bid{doi={10.1214/aop/1024404518}, issn={0091-1798}, mr={1457624}}
\end{barticle}
%
\bptok{imsref}%
% NOT OUTPUTED:
% issn = 0091-1798
% url = http://dx.doi.org/10.1214/aop/1024404518
% number = 3
% coden = APBYAE
% fjournal = The Annals of Probability
\endbibitem

%b19 ###
%b19 #&#
\bibitem{HP98}
%
\begin{barticle}[mr]
\bauthor{\bsnm{H{\"a}ggstr{\"o}m},~\bfnm{Olle}\binits{O.}} \AND
\bauthor{\bsnm{Pemantle},~\bfnm{Robin}\binits{R.}}
(\byear{1998}).
\btitle{First passage percolation and a model for competing spatial growth}.
\bjournal{J. Appl. Probab.}
\bvolume{35}
\bpages{683--692}.
\bid{issn={0021-9002}, mr={1659548}}
\end{barticle}
%
\bptok{imsref}%
% NOT OUTPUTED:
% issn = 0021-9002
% number = 3
% coden = JPRBAM
% fjournal = Journal of Applied Probability
\endbibitem

%b20 ###
%b20 #&#
\bibitem{hueter00}
%
\begin{barticle}[mr]
\bauthor{\bsnm{Hueter},~\bfnm{Irene}\binits{I.}} \AND
\bauthor{\bsnm{Lalley},~\bfnm{Steven~P.}\binits{S.~P.}}
(\byear{2000}).
\btitle{Anisotropic branching random walks on homogeneous trees}.
\bjournal{Probab. Theory Related Fields}
\bvolume{116}
\bpages{57--88}.
\bid{doi={10.1007/PL00008723}, issn={0178-8051}, mr={1736590}}
\end{barticle}
%
\bptok{imsref}%
% NOT OUTPUTED:
% issn = 0178-8051
% url = http://dx.doi.org/10.1007/PL00008723
% number = 1
% coden = PTRFEU
% fjournal = Probability Theory and Related Fields
\endbibitem

%b21 ###
%b21 #&#
\bibitem{KS09}
%
\begin{bincollection}[mr]
\bauthor{\bsnm{Kaimanovich},~\bfnm{Vadim~A.}\binits{V.~A.}} \AND
\bauthor{\bsnm{Sobieczky},~\bfnm{Florian}\binits{F.}}
(\byear{2010}).
\btitle{Stochastic homogenization of horospheric tree products}.
In \bbooktitle{Probabilistic Approach to Geometry}.
\bseries{Adv. Stud. Pure Math.}
\bvolume{57}
\bpages{199--229}.
\bpublisher{Math. Soc. Japan},
\blocation{Tokyo}.
\bid{mr={2648261}}
\end{bincollection}
%
\bptok{imsref}%
\endbibitem

%b22 ###
%b22 #&#
\bibitem{Kersting86}
%
\begin{barticle}[mr]
\bauthor{\bsnm{Kersting},~\bfnm{G.}\binits{G.}}
(\byear{1986}).
\btitle{On recurrence and transience of growth models}.
\bjournal{J. Appl. Probab.}
\bvolume{23}
\bpages{614--625}.
\bid{issn={0021-9002}, mr={0855369}}
\end{barticle}
%
\bptok{imsref}%
% NOT OUTPUTED:
% issn = 0021-9002
% number = 3
% coden = JPRBAM
% fjournal = Journal of Applied Probability
\endbibitem

%b23 ###
%b23 #&#
\bibitem{KL05}
%
\begin{barticle}[mr]
\bauthor{\bsnm{Kordzakhia},~\bfnm{George}\binits{G.}} \AND
\bauthor{\bsnm{Lalley},~\bfnm{Steven~P.}\binits{S.~P.}}
(\byear{2005}).
\btitle{A two-species competition model on {$\mathbb{Z}\sp d$}}.
\bjournal{Stochastic Process. Appl.}
\bvolume{115}
\bpages{781--796}.
\bid{doi={10.1016/j.spa.2004.12.003}, issn={0304-4149}, mr={2132598}}
\end{barticle}
%
\bptok{imsref}%
% NOT OUTPUTED:
% issn = 0304-4149
% url = http://dx.doi.org/10.1016/j.spa.2004.12.003
% number = 5
% coden = STOPB7
% fjournal = Stochastic Processes and their Applications
\endbibitem

%b24 ###
%b24 #&#
\bibitem{Kozma10}
%
\begin{barticle}[mr]
\bauthor{\bsnm{Kozma},~\bfnm{Gady}\binits{G.}}
(\byear{2011}).
\btitle{Percolation on a product of two trees}.
\bjournal{Ann. Probab.}
\bvolume{39}
\bpages{1864--1895}.
\bid{doi={10.1214/10-AOP618}, issn={0091-1798}, mr={2884876}}
\end{barticle}
%
\bptok{imsref}%
% NOT OUTPUTED:
% issn = 0091-1798
% url = http://dx.doi.org/10.1214/10-AOP618
% number = 5
% coden = APBYAE
% fjournal = The Annals of Probability
\endbibitem

%b25 ###
%b25 #&#
\bibitem{LaSe97}
%
\begin{barticle}[mr]
\bauthor{\bsnm{Lalley},~\bfnm{Steven~P.}\binits{S.~P.}} \AND
\bauthor{\bsnm{Sellke},~\bfnm{Tom}\binits{T.}}
(\byear{1997}).
\btitle{Hyperbolic branching {B}rownian motion}.
\bjournal{Probab. Theory Related Fields}
\bvolume{108}
\bpages{171--192}.
\bid{doi={10.1007/s004400050106}, issn={0178-8051}, mr={1452555}}
\end{barticle}
%
\bptok{imsref}%
% NOT OUTPUTED:
% issn = 0178-8051
% url = http://dx.doi.org/10.1007/s004400050106
% number = 2
% coden = PTRFEU
% fjournal = Probability Theory and Related Fields
\endbibitem

%b26 ###
%b26 #&#
\bibitem{Lyons00}
%
\begin{barticle}[mr]
\bauthor{\bsnm{Lyons},~\bfnm{Russell}\binits{R.}}
(\byear{2000}).
\btitle{Phase transitions on nonamenable graphs}.
\bjournal{J. Math. Phys.}
\bvolume{41}
\bpages{1099--1126}.
%statistical physics}.
\bid{doi={10.1063/1.533179}, issn={0022-2488}, mr={1757952}}
\end{barticle}
%
\bptok{imsref}%
% NOT OUTPUTED:
% issn = 0022-2488
% url = http://dx.doi.org/10.1063/1.533179
% number = 3
% coden = JMAPAQ
% fjournal = Journal of Mathematical Physics
\endbibitem

%b27 ###
%b27 #&#
\bibitem{woess}
%
\begin{bbook}[mr]
\bauthor{\bsnm{Woess},~\bfnm{Wolfgang}\binits{W.}}
(\byear{2000}).
\btitle{Random Walks on Infinite Graphs and Groups}.
\bseries{Cambridge Tracts in Mathematics}
\bvolume{138}.
\bpublisher{Cambridge Univ. Press},
\blocation{Cambridge}.
\bid{doi={10.1017/CBO9780511470967}, mr={1743100}}
\end{bbook}
%
\bptok{imsref}%
% NOT OUTPUTED:
% isbn = 0-521-55292-3
% url = http://dx.doi.org/10.1017/CBO9780511470967
% fpage = xii+334
\endbibitem

%b28 ###
%b28 #&#
\bibitem{Zucca11}
%
\begin{barticle}[mr]
\bauthor{\bsnm{Zucca},~\bfnm{Fabio}\binits{F.}}
(\byear{2011}).
\btitle{Survival, extinction and approximation of discrete-time
branching random walks}.
\bjournal{J. Stat. Phys.}
\bvolume{142}
\bpages{726--753}.
\bid{doi={10.1007/s10955-011-0134-x}, issn={0022-4715}, mr={2773785}}
\end{barticle}
%
\bptok{imsref}%
% NOT OUTPUTED:
% issn = 0022-4715
% url = http://dx.doi.org/10.1007/s10955-011-0134-x
% number = 4
% fjournal = Journal of Statistical Physics
\endbibitem

\end{thebibliography}
\end{document}